\documentclass[12pt]{article}

\usepackage{amssymb,amsmath}

\newcommand{\qed}{{\hfill$\blacksquare$}}
\newcommand{\CC}{\mathbb C}

\newcommand{\ZZ}{\mathbb Z}
\newcommand{\PP}{\mathbb P}

\newcommand{\Ma}{\mathsf{M}}
\newcommand{\ma}{\mathsf{m}}
\newcommand{\EC}{\overline{\mathcal{E}}}
\newcommand{\ECF}{\mathcal{E}}

\newcommand{\tP}{\widetilde{P}}

\newcommand{\nr}{{$\sharp$ }}

\newcommand{\gl}{{\lambda }}
\newcommand{\eps}{{\epsilon }}
\newcommand{\gx}{{\xi }}

\newcommand{\ga}{{\alpha }}

\newcommand{\CQ}{{\mathsf{Q} }}
\newcommand{\cq}{{\mathsf{q} }}

\begin{document}

\title{Mahler Measure Variations, Eisenstein Series and Instanton Expansions.}

\author{Jan Stienstra\footnote{e-mail: stien{@}math.uu.nl}\\
\small{Mathematisch Instituut, Universiteit Utrecht, the Netherlands}}
\date{}

\maketitle

\begin{abstract}
This paper points at an intriguing inverse function relation with
on the one hand the coefficients of the Eisenstein series in Rodriguez Villegas' paper on ``Modular Mahler Measures'' and on the other hand the instanton numbers in papers on ``Non-Critical Strings'' by Klemm-Mayr-Vafa and Lerche-Mayr-Warner. In a companion paper Mahler measures will be related to dimer models. Thus, if it is not just a lucky coincidence in a few examples, our inverse function relation may be an incarnation of the duality between string models and dimer models proposed by Okounkov, Reshetikhin and Vafa in their paper ``Quantum Calabi-Yau and Classical Crystals''.
\end{abstract}

\section*{Introduction}

The \emph{logarithmic Mahler measure} $\ma (F)$ and the \emph{Mahler measure} $\Ma (F)$ of a Laurent polynomial $F(x,y)$ with complex coefficients are:
\begin{eqnarray}\label{eq:logmahler}
\ma (F)&:=&\frac{1}{(2\pi i)^2}\oint\!\!\oint_{|x|=|y|=1}
\log |F(x,y)|\,
\frac{dx}{x}\frac{dy}{y}\;,\\
\label{eq:mahler}
\Ma (F)&:=&\exp (\ma (F))\;.
\end{eqnarray}
Wondering how the Mahler measure depends on the coefficients of the Laurent polynomial, and more in particular on the constant term,  we will look at the Mahler measure $\Ma(F_t)$ of Laurent polynomials
$$
F_t(x,y):=t-(xy)^{-1}\tP(x,y,1)
$$ 
where $\tP(X,Y,Z)$ is a fixed weighted homogeneous polynomial and
$t$ is a complex variable. Geometrically speaking, we
look at the Mahler measure of the fibres of the map
$\PP^2\longrightarrow \PP^1$,
$(X,Y,Z)\mapsto[\tP(X,Y,Z):(XYZ)]$
on the (weighted) projective plane $\PP^2$. In the examples these fibres are elliptic curves. The variation in the Mahler measure is then related to the variation in the periods, or rather the period ratio, of these elliptic curves 
(see \S\ref{Section 1}).
This is what Rodriguez Villegas was investigating when he noticed that some of his formulas looked similar to formulas in the literature on Mirror Symmetry; see the final remarks in \cite{RV}.
The modest goal of the present paper is to point out that actual instanton expansions matching his examples do occur in the physics literature \cite{LMW, KMV, KZ, CKYZ}. The relation is, however, more intriguing and less straightforward than one might expect from the remarks in \cite{RV}:

\

\noindent
\textbf{Observation:}\footnote{I call this an observation and not a theorem, since it can be `proved' 
by just pointing out how pieces of existing literature fit together
(which, however,
as far as I know, were never presented together in this form).}
\textit{
The expansions $1-\sum b_n\frac{n^2 \cq^n}{1-\cq^n}$ of the weight $3$ modular forms in
Examples $1,\; 2,\; 3$ of \cite{RV} are related to 
the instanton expansions 
\mbox{$\kappa(1-\sum a_n\frac{n^3 \CQ^n}{1-\CQ^n})$} in
Examples $E_6,\;E_5,\;E_8$, respectively, in \cite{KMV, LMW}, via the inverse function relations}
\begin{equation}
\label{eq:Qqintro}
-\CQ=\cq\prod_{n\geq 1} (1-\cq^n)^{nb_n}\:,\qquad
-\cq=\CQ\prod_{n\geq 1} (1-\CQ^n)^{n^2a_n}\:.
\end{equation}
\textit{The constant $\kappa$ is $-3, -4, -1$ for Examples $E_6,\;E_5,\;E_8$, respectively.}

\textit{
The Mahler measure of the elliptic curve with modular parameter $\cq$ is $|\CQ|^{-1/\nu}$ with $\nu=3,\, 2,\, 6$ for Examples $E_6,\;E_5,\;E_8$, respectively.}
\qed

\

That Examples $1,\; 2,\; 3$ of \cite{RV} line up with Examples $E_6,\;E_5,\;E_8$ in \cite{KMV, LMW} is immediately obvious from the formulas for the elliptic pencils and their Picard-Fuchs equations which are explicitly given in \cite{KMV, LMW, RV}. Looking at how the various authors get from the Picard-Fuchs equations to their expansions  directly 
leads to Relation (\ref{eq:Qqintro}); for more details see Section \ref{Section 1}.

In Section \ref{section examples}
we present six such examples including the above ones and the $E_7$ case of \cite{KMV, LMW}.
In each case the computer output suggests that the numbers
$a_n$ and $b_n$ are integers. Except for the $b_n$ in four  examples we have no proof for this suggested integrality.

\

\noindent
\textbf{Discussion and new directions.} While the second product in (\ref{eq:Qqintro}) summarizes some of results (instanton expansions) in the String Theory literature of the 1990's, Mahler measure and an analogue of the first product in (\ref{eq:Qqintro}) appear in the String related literature of the 2000's. Theorem 3.5 of \cite{KOS}, for instance, gives a formula for the partition function of a dimer model which is exactly the same as the defining formula 
(\ref{eq:mahler}) for the Mahler measure of the characteristic polynomial of that model. In \cite{S2} we show that the polynomials of Examples \nr 3, \nr 4 and \nr 6 of the present paper are characteristic polynomials of dimer models.
The partition functions of these dimer models can therefore be given as
the multiplicative inverses of the first products in (\ref{eq:Qqintro}). For Example \nr 3, for instance, this product is
$$
\cq\prod_{n\geq 1}(1-\cq^n)^{9n\chi_{-3}(n)}\,,\qquad
\chi_{-3}(n)=0,1,-1\quad\textrm{for}\quad n\equiv 0, 1, 2 \bmod 3.
$$
There is here some similarity with McMahon's function
$$
M(q):=\prod_{n\geq 1}(1-\cq^n)^{-n}
$$
which appears in partition functions in \cite{O, ORV}. The fact that McMahon's function is the generating function for plane partitions (i.e. with 3 dimensional Young diagrams) led to unexpected new directions in String Theory research: \emph{the melting crystal and random surface pictures} \cite{KOS, O, ORV}.

\

\noindent
\textbf{Acknowledgements.} I want to thank Don Zagier for bringing
\cite{RV} to my attention (somewhen around 1998) and for keeping me well informed about his work in progress \cite{Z}. I thank
Bogdan Florea for discussing with me in November 2001 some parts of his work, thereby also providing references to \cite{KMV, LMW}. 
I want to thank Noriko Yui for the invitation to the Banff Workshop on Mirror Symmetry with its stimulating atmosphere to report on this work.

\section{Details of the Observation} \label{Section 1}

We investigate Laurent polynomials of the form  
$$
F_t(x,y):=t-(xy)^{-1}\tP(x,y,1)
$$ 
with $t$ a complex parameter and $\tP(X,Y,Z)$ as in Table 
\ref{table 1}.
Writing $\gx=t^{-1}$ and assuming 
$\max_{|x|=|y|=1} |\tP(x,y,1)||\gx|<1$
one then finds 
\begin{equation}
\label{eq:mahlerQ}
\Ma(F_t)=|Q|^{-1}
\end{equation}
with
\begin{equation}
\label{eq:period1}
\begin{array}{rcl}
Q&:=&\displaystyle{\beta\,
\exp\left(-\frac{1}{(2\pi i)^2}\oint\!\!\oint_{|x|=|y|=1}
\log (t-\frac{\tP(x,y,1)}{xy})\,
\frac{dx}{x}\frac{dy}{y}\right)}
\\[1.5ex]
&=&\displaystyle{\beta\,\gx \exp\left(\sum_{n=1}^\infty
\frac{\gx^n}{n}
\left(\frac{1}{(2\pi i)^2}\oint\!\!\oint_{|x|=|y|=1}
\frac{\tP(x,y,1)^n}{(xy)^n}\,\frac{dx}{x}\frac{dy}{y}\right)\right)}
\end{array}
\end{equation}
and $\beta$ an as yet undetermined complex number of absolute value $1$.
The number $\beta$ is, of course, irrelevant for (\ref{eq:mahlerQ}), but
for the correct match with \cite{LMW, KMV} it will be necessary to take here an appropriate $\gx$-independent constant.

Thus
\begin{equation}
\label{eq:Qxi}
Q\,=\,\beta\,\gx\,\exp\left(\sum_{n=1}^\infty c_n\frac{\gx^n}{n}\right)
\end{equation}
with
\begin{equation}
\label{eq:cn}
c_n =
\textrm{coefficient of } X^nY^nZ^n \textrm{ in }\: \tP(X,Y,Z)^n.
\end{equation}
We refer to Table \ref{table 2} for concrete examples.
Now let $\EC_t$ denote the closure of
$$ 
\ECF_t=\{ (x,y)\in \CC^*\times\CC^*\:|\:\tP(x,y,1)=txy\:\}
$$ 
in the (appropriately weighted) projective plane.
Restricting to values of $t=\gx^{-1}$ for which the curve $\EC_t$
is non-singular, 
differentiating (\ref{eq:period1}) and using the residue theorem one finds
\begin{equation}
\label{eq:period2}
\gx\frac{d}{d\gx}\log Q
\;=\;\frac{1}{(2\pi i)^2}\oint\!\!\oint_{|x|=|y|=1}
\frac{t\,dx\,dy}{t\,xy- \tP(x,y,1)}
\;=\;
\frac{1}{2\pi i}\oint_{\gamma}\omega_t
\end{equation}
where $\gamma$ is a closed loop and $\omega_t$ is a holomorphic 1-form on the curve $\EC_t$,
$$
\omega_t:=\frac{tdx}{tx- \frac{\partial}{\partial y}\tP(x,y,1)}
$$
on the coordinate patch where the denominator does not vanish.

Thus $\gx\frac{d}{d\gx}\log Q$ is a period of $\omega_t$.
The periods of $\omega_t$ are functions of $\gx=t^{-1}$ which satisfy a second order linear differential equation: the \emph{Picard-Fuchs equation} of the elliptic pencil. In our examples there is only one 
solution $g_1(\gx)$ of the Picard-Fuchs equation which is holomorphic in a neighborhood of $\gx=0$ and has $g_1(0)=1$.
Differentiating (\ref{eq:Qxi}) shows that this solution is
\begin{equation}
\label{eq:period3}
g_1(\gx)\;=\;\gx\frac{d}{d\gx}\log Q\;=\;\sum_{n=0}^\infty c_n\gx^n.
\end{equation}
A second solution has the form
$$
g_2(\gx)\;=\;g_1(\gx)\log\gx\:+\: h(\gx),
$$
where $h(\gx)$ is a holomorphic function of $\gx$ which vanishes at $\gx=0$. The classical procedure in analyzing elliptic pencils is to consider the period ratio and to set
\begin{equation}
\label{eq:qxi}
q\: :=\:\exp\left(\frac{g_2(\gx)}{g_1(\gx)}\right)
\:=\: \gx\exp\left(\frac{h(\gx)}{g_1(\gx)}\right).
\end{equation}
Formulas (\ref{eq:Qxi}) and (\ref{eq:qxi}) give $Q$ and $q$ as functions of $\gx$, while
$$
q=\gx + O(\gx^2)\,,\qquad Q=\beta\gx+ O(\gx^2).
$$
One can eliminate $\xi$ and directly express $Q$ as a function of $q$ and vice versa. Although, one could now express this inverse function relationship through product expansions as in (\ref{eq:Qqintro}), it is better to wait: in the examples there is a positive integer $\nu$ (see Table \ref{table 2}) such that $Q^\nu$ is a function of $q^\nu$ and vice versa. It is this latter inverse function relationship which is used in (\ref{eq:Qqintro}) with
\begin{equation}\label{eq:QQqq}
\CQ:=Q^\nu\,,\qquad \cq:=q^\nu\,.
\end{equation}
Here are more details about the role of $\nu$.
One solution of the Picard-Fuchs differential equation is:
\begin{equation}\label{eq:generating}
g_1=\sum_{m=0}^\infty u_m\psi^m=\sum_{n=0}^\infty c_n\gx^n.
\end{equation}
with $\psi=\gx^\nu$, $u_m$ and $c_n$ as in Table \ref{table 2}. The 
numbers $u_m$ and $(A,B,\gl)$ in Table \ref{table 2}
satisfy the recurrence relation
\begin{equation}\label{eq:recurrence}
(m+1)^2u_{m+1}-(Am^2+Am+\gl)u_m+Bm^2u_{m-1}= 0\qquad (m\geq 0)\,.
\end{equation}
This can easily be checked by hand for the 
hypergeometric cases \nr 1--\nr 4. The recurrence in case \nr 5 is due to Ap\'ery; the numbers $u_m$ appear in his famous irrationality proof for $\zeta(2)$. The recurrence in case \nr 6 is due to Cusick.
Proofs for the recurrences in cases \nr 5 and \nr 6 (which in fact first derive the Picard-Fuchs equation) can be found in, for instance, \cite{SB} \S 11. The $(A,B,\gl)$-notation was introduced by
Zagier \cite{Z}, who made an extensive search for all such triples for which all $u_m$ in the recurrence (\ref{eq:recurrence}) are integers.

The recurrence relation (\ref{eq:recurrence}) is equivalent with the differential equation
\begin{equation}\label{eq:differential}
\psi\frac{d}{d\psi}\left((1-A\psi+B\psi^2)\psi\frac{d}{d\psi} g\right) +\psi(-\gl+B\psi)g=0
\end{equation}
for the generating function $g_1$ in (\ref{eq:generating}).
A second solution of the differential equation (\ref{eq:differential})
can be obtained by the method of Frobenius. For this we
make the Ansatz
$$
g(\psi)=\sum_{m=0}^\infty u_m(\eps)\psi^{m+\eps}\,,\qquad u_0(\eps)=1.
$$ 
This is a solution of (\ref{eq:differential}) if and only if $\eps^2=0$
and for all $m\geq 0$
\begin{equation}\label{eq:epsrecur}
(m+\eps)^2u_m(\eps)=(A(m+\eps)^2-A(m+\eps)+\gl)u_{m-1}(\eps)
-B(m-1+\eps)^2u_{m-2}(\eps).
\end{equation}
We set $u_0(\eps)=1$, solve the recurrrence and then
expand 
\begin{equation}\label{eq:gpsi}
g(\psi)=g_1(\psi)+g_2(\psi)\eps.
\end{equation}
Then $g_1(\psi)$ and $g_2(\psi)$ are solutions of (\ref{eq:differential}),
$g_1(\psi)=\sum_{m=0}^\infty u_m\psi^m$ as before and
$g_2(\psi)=g_1(\psi)\log\psi+h(\psi)$
where $h(\psi)$ is a power series in $\psi$ with $h(0)=0$.
Now set, with $\nu$ as in Table \ref{table 2},
\begin{equation}\label{eq:q}
q:=\exp\left(\frac{g_2(\psi)}{\nu\,g_1(\psi)}\right)=
\gx\exp\left(\frac{h(\gx^\nu)}{\nu\,g_1(\gx^\nu)}\right)\,.
\end{equation}
This definition of $q$ agrees with (\ref{eq:qxi}), but note that
$\cq:=q^\nu$ is actually a function of $\psi=\gx^\nu$ and
$\cq=\psi+O(\psi^2)$.

Now recall (\ref{eq:period3}) and (\ref{eq:generating}):
\begin{equation}\label{eq:dQdpsi}
\psi\frac{d}{d\psi}\log Q^\nu\;=\;
\gx\frac{d}{d\gx}\log Q\;=\;\sum_{n=0}^\infty c_n\gx^n\;=\;\sum_{m=0}^\infty u_m\psi^m\,.
\end{equation}
This shows that $\CQ:=Q^\nu$ is actually a function of $\psi$ and that
$\CQ=\ga\psi+O(\psi^2)$ for some constant $\ga$.
One can eliminate $\psi$ and directly express $\CQ$ as a function of $\cq$ and vice versa.
This is the inverse function relationship used in (\ref{eq:Qqintro}),
with $\ga=-1$.

By comparing the above formulas for $\CQ$ and for $\cq$
with the formulas in $\textrm{n}^\circ$13 of
\cite{RV} one sees that our 
$\cq\frac{d}{d\cq}\log \CQ=\left(\psi\frac{d}{d\psi}\log \CQ\right)\left(\cq\frac{d}{d\cq}\log \psi\right)$
equals Rodriguez Villegas' $e(\tau)$, which in the examples has an 
expansion like an Eisenstein series of weight $3$:
\begin{equation}\label{eq:dQdq}
\cq\frac{d}{d\cq}\log \CQ=1-\sum_{n\geq 1} b_n\frac{n^2 \cq^n}{1-\cq^n}\,.
\end{equation}
Note that (\ref{eq:period3}) determines $\CQ$ only up to a multiplicative constant. This constant is irrelevant for the expansion (\ref{eq:dQdq}). 
On the other hand, for \cite{LMW} fixing this constant is very relevant.
It is called fixing the B-field. 
Equations (4.10) and (4.13) of \cite{LMW} essentially boil down to
\begin{equation}\label{eq:dqdQ}
\CQ\frac{d}{d\CQ}\log \cq\:=\:1-\sum_{n\geq 1} a_n\frac{n^3 \CQ^n}{1-\CQ^n}\,;
\end{equation}
to really get the expansions of \cite{LMW} Equation (4.10) the right hand side of (\ref{eq:dqdQ}) has to be multiplied by an appropriate constant
(see the numerical results in the next section).

With (\ref{eq:dQdq}) and (\ref{eq:dqdQ}) one can express $\CQ$ as a function of $\cq$ and vice versa, both up to a multiplicative constant. These constants are fixed by 
$$
\CQ=\ga\psi+O(\psi^2)\,,\qquad\cq=\psi+O(\psi^2)
$$
Thus one arrives
at the following functional relation between $\CQ$ and $\cq$:
\begin{equation}
\label{eq:Qq}
\ga^{-1} \CQ=\cq\prod_{n\geq 1} (1-\cq^n)^{nb_n}\:,\qquad
\ga \cq=\CQ\prod_{n\geq 1} (1-\CQ^n)^{n^2a_n},
\end{equation}
Incidentally, the relation between the $\ga$ here and the $\beta$ in
(\ref{eq:period1}) is:  $\ga=\beta^\nu$.
In the examples $\ga=\pm 1$ and all $a_n$ 
and $b_n$ turn out to be integers. Except for some of the $b_n$ we have no proof for this integrality.

\section{Numerical results for six examples.}\label{section examples}
Our examples deal with pencils of elliptic curves in a (weighted) projective plane $\PP$ described by the polynomials in 
Table \ref{table 1}:
\begin{equation}\label{eq:compact ell. curve2}
\EC_t:=\{(X,Y,Z)\in\PP\:|\: \tP(X,Y,Z)=tXYZ\}.
\end{equation}

Beauville \cite{Beau} showed that there are exactly six semi-stable families of elliptic curves over $\PP^1$ with four singular fibres.
His list consists of cases \mbox{\nr 3 - \nr 6} in our Table \ref{table 1}
plus the polynomials
 $ (X+Y)((X+Y+Z)^2-4YZ) $ and \mbox{$ X^3+Y^3+Z^3 $}. We do not explicitly discuss the latter two, since the results coincide with those of \nr 4 resp. \nr 3.

Table \ref{table 2} shows that cases \mbox{\nr 1 - \nr 4} are of hypergeometric origin; more precisely they result by a substution
$t\mapsto t^6,t^4,t^3,t^2$, respectively, from elliptic pencils with exactly three singular fibres. Families of elliptic curves with exactly three singular fibres have been classified in \cite{SH}. From the tables in \cite{SH} one can see that our examples come from the singular fibre configurations
$\sharp 1:(II^*\;I_1\;I_1)\,,\;
\sharp 2:(III^*\; I_2\; I_1)\,,\;
\sharp 3:(IV^*\;I_1\;I_3)\,,$
\mbox{$\sharp 4:(I^*_1\;I_1\;I_4)$.}
Kodaira's notation for the singular fibres corresponds as follows with Dynkin diagrams:
$II^*= E_8$,
$III^*= E_7$,
$IV^*= E_6$,
$I^*_1= E_5=D_5$.
This agrees with the computations hereafter according to which examples \mbox{\nr 1 --\nr 4} give cases $E_8,E_7,E_6,E_5$ in \cite{KMV,LMW}.

\begin{table}
\begin{center}
\begin{tabular}{|l|cl|}\hline 
\rule{0mm}{4mm}
& &\hspace{3em} $\tP(X,Y,Z)$
\\
\hline
\nr  1&\rule{0mm}{4mm}\hspace{3em}&
$ X^2+Y^3+Z^6 $
\\
\nr  2& &
$ X^2+Y^4+Z^4 $
\\
\nr  3& &
$ X^2Y+Y^2Z+Z^2X $ 
\\
\nr  4& &
$ (X+Y)(XY+Z^2) $
\\
\nr  5& &
$ (X+Y+Z)(X+Z)(Y+Z) $
\\
\nr  6& &
$ (X+Y)(Y+Z)(Z+X) $
\\
\hline
\end{tabular}
\caption{}\label{table 1}
\end{center}
\end{table}

\begin{table}
\begin{center}
\begin{tabular}{|l|cc|cr|cc|cc|}\hline 
\rule{0mm}{4mm}
& \hspace{3em}& $u_m$&\hspace{3em}&$c_n=u_m\hspace{1em}$& \hspace{3em}& $(A,B,\gl)$&\hspace{3em}&$\nu$
\\
\hline
\nr  1&\rule{0mm}{4mm}\hspace{5em}&
$\displaystyle{ \frac{(6m)!}{m!(2m)!(3m)!} }$& & if $6|n\,$, 
$\; m=\frac{n}{6}$ & &(432,\,0,\,60)& &$6$
\\[2ex]
\nr  2& &
$\displaystyle{ \frac{(4m)!}{m!^2(2m)!} }$& & if $4|n\,$, 
$\; m=\frac{n}{4}$ & &(64,\,0,\,12)& &$4$
\\[2ex]
\nr  3& &
$\displaystyle{ \frac{(3m)!}{m!^3} }$& & if $3|n\,$, 
$\; m=\frac{n}{3}$ & &(27,\,0,\,6)& &$3$
\\[2ex]
\nr  4& &
$\displaystyle{ \frac{(2m)!^2}{m!^4} }$& & if $2|n\,$, 
$\; m=\frac{n}{2}$ & &(16,\,0,\,4)& &$2$
\\[3ex]
\nr  5& &
$
\displaystyle{\sum_{k=0}^m} {\scriptscriptstyle{
\left(\begin{array}{c}m\\ k\end{array}\right)}}^{\scriptstyle{2}}
\scriptscriptstyle{\left(\begin{array}{c}m+k\\ k\end{array}\right)
}$& & $m=n$ & &(11,\,-1,\,3)& &$1$
\\[4ex]
\nr  6& &
$
\displaystyle{\sum_{k=0}^m} {\scriptscriptstyle{
\left(\begin{array}{c}m\\ k\end{array}\right)}}^{\scriptstyle{3}}$& & $m=n$
& &(7,\,-8,\,2)& &$1$
\\
\hline
\end{tabular}
\caption{$c_n=0$ if $n$ does not satisfy the condition in the 
third column }\label{table 2}
\end{center}
\end{table}

We now come to the actual calculations (using PARI). 
From the input in the form of the triple $(A,B,\gl)$ we compute $g_1$ and $\cq=q^\nu$ as functions of $\psi$ using Formulas (\ref{eq:epsrecur})-(\ref{eq:q}). Inverting the latter functional relation we get $\psi$ and $g_1$ as functions of $\cq$. The numbers $b_n$ are then computed from the expansion
(cf. (\ref{eq:dQdq}), (\ref{eq:dQdpsi}))
\begin{equation}\label{eq:eisq}
g_1(\cq)\,\frac{\cq}{\psi(\cq)}\:\frac{d\psi(\cq)}{d\cq}\;=\;
1-\sum_{n\geq 1} b_n\frac{n^2 \cq^n}{1-\cq^n}\;=\;
1-\sum_{n\geq 1}\left(\sum_{k|n} k^2b_k\right)\cq^n\,.
\end{equation}
Having the numbers $b_n$ we compute $\CQ$ and the numbers $a_n$ from
(\ref{eq:Qq}):
\begin{equation}\label{eq:Qqrel}
\ga^{-1} \CQ=\cq\prod_{n\geq 1} (1-\cq^n)^{nb_n}\:,\qquad
\ga \cq=\CQ\prod_{n\geq 1} (1-\CQ^n)^{n^2a_n}\,.
\end{equation}
Typically $\ga=-1$ leads to the numbers which are reported 
in 
\cite{LMW} Equations (4.10); for the conifold case, though, $\ga=1$,
as remarked in loc. cit..

We report the results of the calculations in tables for the numbers $a_n$ and $b_n$ and as $\eta$-products or Eisenstein series for
$\psi(\cq)$, $g_1(\psi(\cq))$ and $\cq\frac{d}{d\cq}\log\CQ$. The elliptic pencils and Picard-Fuchs equations of Examples \nr 3-\nr 6 are well-known and well documented in the literature \cite{Z, RV}; also \nr 1 appears in \cite{RV} as example 3. So we have taken the expressions for $\psi(\cq)$ and $g_1(\psi(\cq))$ directly from \cite{Z, RV}.
Zagier \cite{Z} gives an argument for ``recognizing modularity'' which actually proves that the expressions are correct.
However, as an extra security check we have checked that the $\cq$-expansions of these expressions coincide with our PARI calculations to order at least $100$.

\subsection*{\nr 1: $(A,B,\gl)=(432,0,60)\,,\;\ga=-1$.}

This corresponds to the $E_8$ example in \cite{LMW,KMV} and Example 3 in \cite{RV}. 
$$
\begin{array}{rrrrr}
n&\qquad&a_n&\qquad&b_n\\
&&&&\\
1&& 252&& 252\\
2&&-9252&&-13374\\
3&& 848628&& 1253124\\
4&&-114265008&&-151978752\\
5&& 18958064400&& 21255487740\\
6&&-3589587111852&&-3255937602498\\
7&& 744530011302420&& 531216722607876\\
8&&-165076694998001856&&-90773367805541376\\
9&& 38512679141944848024&& 16069733941012586748
\end{array}
$$
The numbers of $a_n$ are the instanton numbers of \cite{LMW,KMV} Example $E_8$. They also appear in \cite{CKYZ} Table 7
column $X_6(1123)$. More precisely the expansion for the $E_8$ case
in \cite{LMW} Eq. (4.10) is:
$$
-1\cdot\left(1-\sum_{n\geq 1} a_n\frac{n^3 \CQ^n}{1-\CQ^n}\right)
$$
Rodriguez Villegas \cite{RV} Example 3 gives the following expressions
for, in our notation, $\psi(\cq)$, $g_1(\psi(\cq))$ and 
$\cq\frac{d}{d\cq}\log\CQ$:
\begin{eqnarray*}
\psi(\cq)&=&\frac{1}{864}\left(
1-E_6(\cq)E_4(\cq)^{-3/2}\right)\\
g_1(\psi(\cq))&=&E_4(\cq)^{1/4}\\
\cq\frac{d}{d\cq}\log\CQ&=&1-\sum_{n\geq 1} b_n \frac{n^2\cq^n}{1-\cq^n}=\frac{1}{2}
E_4(\cq)^{3/4}\left(1+E_6(\cq)E_4(\cq)^{-3/2}\right)
\end{eqnarray*}
where $E_4(\cq)$ and $E_6(\cq)$ are the Eisenstein series
$$
E_4(\cq):=1+240\sum_{n=1}^\infty\sum_{d|n}d^3\cq^n
\,,\qquad
E_6(\cq):=1-504\sum_{n=1}^\infty\sum_{d|n}d^5\cq^n\,.
$$

\subsection*{\nr 2: $(A,B,\gl)=(64,0,12)\,,\;\ga=-1$.}
This corresponds to the $E_7$ example in \cite{LMW,KMV}, but has no explicitly worked out counterpart in \cite{RV}.
$$
\begin{array}{rrrrr}
n&\qquad&a_n&\qquad&b_n\\
&&&&\\
1&& 28&& 28\\
2&&-136&&-134\\
3&& 1620&& 996\\
4&&-29216&&-10720\\
5&& 651920&& 139292\\
6&&-16627608&&-2019450\\
7&& 465215604&& 31545316\\
8&&-13927814272&&-520076672\\
9&& 439084931544&& 8930941980
\end{array}
$$ 
The numbers of $2a_n$ are the instanton numbers of \cite{LMW,KMV} Example $E_7$. They also appear in \cite{CKYZ} Table 7
column $X_4(1112)$. More precisely the expansion for the $E_7$ case
in \cite{LMW} Eq. (4.10) is:
$$
-2\cdot\left(1-\sum_{n\geq 1} a_n\frac{n^3 \CQ^n}{1-\CQ^n}\right)
$$
Here we can only offer (without reference to proofs in the literature) expressions for $\psi(\cq)$, $g_1(\psi(\cq))$ and 
$\cq\frac{d}{d\cq}\log\CQ$, which are correct to at least $O(\cq^{100})$:
\begin{eqnarray*}
\psi(\cq)^{-1}&=&64+\frac{\eta(\cq)^{24}}{\eta(\cq^2)^{24}}\\
g_1(\psi(\cq)&=&(2E_2(\cq^2)-E_2(\cq))^{\frac{1}{2}}\\
\cq\frac{d}{d\cq}\log\CQ&=&
1-\sum_{n\geq 1} b_n \frac{n^2\cq^n}{1-\cq^n}=
(2E_2(\cq^2)-E_2(\cq))^{\frac{3}{2}} \left(1+64\frac{\eta(\cq^2)^{24}}{\eta(\cq)^{24}}\right)^{-1}
\end{eqnarray*}
where 
$$
E_2(\cq):=1-24\sum_{n=1}^\infty\sum_{d|n}d \cq^n
\,,\qquad
\eta(\cq):=\cq^{\frac{1}{24}}\prod_{n\geq 1}(1-\cq^n)\,.
$$

\subsection*{\nr 3: $(A,B,\gl)=(27,0,6)\,,\;\ga=-1$.}
This corresponds to the $E_6$ example in \cite{LMW,KMV} and Example 1 in \cite{RV}. 
$$
\begin{array}{rrrrr}
n&\qquad&a_n&\qquad&b_n\\
&&&&\\
1&& 9&&9\\
2&&-18&&-9\\
3&& 81&&0\\
4&&-576&&9\\
5&& 5085&&-9\\
6&&-51192&&0\\
7&& 565362&&9\\
8&&-6684480&&-9\\
9&& 83246697&&0\\
\end{array}
$$ 
The numbers of $3a_n$ are the instanton numbers of \cite{LMW,KMV} Example $E_6$. The numbers $\frac{1}{3}a_n$
appear in \cite{CKYZ} Table 1 and in \cite{KZ} Fig. 2.
More precisely the expansion for the $E_6$ case
in \cite{LMW} Eq. (4.10) is:
$$
-3\cdot\left(1-\sum_{n\geq 1} a_n\frac{n^3 \CQ^n}{1-\CQ^n}\right)
$$
Zagier \cite{Z} and Rodriguez Villegas \cite{RV} give the following expressions for, in our notation, $\psi(\cq)$, $g_1(\psi(\cq))$ and 
$\cq\frac{d}{d\cq}\log\CQ$:
\begin{eqnarray*}
\psi(\cq)^{-1}&=&27+\frac{\eta(\cq)^{12}}{\eta(\cq^3)^{12}}\\
g_1(\psi(\cq))&=&E_{1,\chi_{-3}}(\cq)\,:=
\,1+6\sum_{n\geq 1} \chi_{-3}(n) \frac{\cq^n}{1-\cq^n}\\
\cq\frac{d}{d\cq}\log\CQ&=&1-\sum_{n\geq 1} b_n \frac{n^2\cq^n}{1-\cq^n}=
E_{3,\chi_{-3}}(\cq)\,:=\,1-9\sum_{n\geq 1} \chi_{-3}(n) \frac{n^2\cq^n}{1-\cq^n}
\end{eqnarray*} 
where $\chi_{-3}$ denotes the quadratic Dirichlet character 
$\left(\frac{-3}{.}\right)$, i.e.
$$
\chi_{-3}(n)=\textstyle{\left(\frac{-3}{n}\right)}=0,\,1,\,-1\qquad
\textrm{if}\quad n\equiv 0,\,1,\,2\bmod 3.
$$ 
According to \cite{ABYZ} one also has
$
E_{3,\chi_{-3}}(\cq)\,=\,\eta(\cq)^9\eta(\cq^3)^{-3}\,.
$

\subsection*{\nr 4: $(A,B,\gl)=(16,0,4)\,,\;\ga=-1$.}
This corresponds to the $E_5$ example in \cite{LMW,KMV} and Example 2 in \cite{RV}. 
$$
\begin{array}{rrrrr}
n&\qquad&a_n&\qquad&b_n\\
&&&&\\
1&& 4&&4\\
2&&-5&&0\\
3&& 12&&-4\\
4&&-48&&0\\
5&& 240&&4\\
6&&-1359&&0\\
7&& 8428&&-4\\
8&&-56000&&0\\
9&& 392040&&4\\
\end{array}
$$  The numbers of $4a_n$ are the instanton numbers of \cite{LMW,KMV} Example $E_5$. More precisely the expansion for the $E_5$ case
in \cite{LMW} Eq. (4.10) is:
$$
-4\cdot\left(1-\sum_{n\geq 1} a_n\frac{n^3 \CQ^n}{1-\CQ^n}\right)
$$
Zagier \cite{Z} and Rodriguez Villegas \cite{RV} give the following expressions for, in our notation, $\psi(\cq)$, $g_1(\psi(\cq))$ and 
$\cq\frac{d}{d\cq}\log\CQ$:
\begin{eqnarray*}
\psi(\cq)&=&\frac{\eta(\cq)^8\eta(\cq^4)^{16}}{\eta(\cq^2)^{24}}\\
g_1(\psi(\cq))&=&E_{1,\chi_{-4}}(\cq)\,:=
\,1+4\sum_{n\geq 1} \chi_{-4}(n) \frac{\cq^n}{1-\cq^n}\\
\cq\frac{d}{d\cq}\log\CQ&=&1-\sum_{n\geq 1} b_n \frac{n^2\cq^n}{1-\cq^n}=
E_{3,\chi_{-4}}(\cq)\,:=\,1-4\sum_{n\geq 1} \chi_{-4}(n) \frac{n^2\cq^n}{1-\cq^n}
\end{eqnarray*} 
where $\chi_{-4}$ denotes the quadratic Dirichlet character 
$\left(\frac{-4}{.}\right)$, i.e.
$$
\chi_{-4}(n)=\textstyle{\left(\frac{-4}{n}\right)}=0,\,1,\,0,\,-1\qquad
\textrm{if}\quad n\equiv 0,\,1,\,2,\,3\bmod 4.
$$ 
According to \cite{ABYZ} one also has
$
E_{3,\chi_{-4}}(\cq)\,=\,\eta(\cq)^4\eta(\cq^2)^6\eta(\cq^4)^{-4}\,.
$

\

\textbf{Remark:} Taking $(A,B,\gl)=(16,0,4)\,,\;\ga=1$
leads to the instanton numbers of \cite{LMW} Example Conifold.

\subsection*{\nr 5: $(A,B,\gl)=(11,-1,3)\,,\;\ga=-1$.}
This example does not appear in \cite{LMW,KMV} or \cite{RV}.
$$
\begin{array}{rrrrr}
n&\qquad&a_n&\qquad&b_n\\
&&&&\\
1&& 2&&2\\
2&&-2&&1\\
3&& 3&&-1\\
4&&-8&&-2\\
5&& 27&&0\\
6&&-102&&2\\
7&& 420&&1\\
8&&-1856&&-1\\
9&& 8649&&-2\\
\end{array}
$$
Zagier \cite{Z} gives for $\psi(\cq)$ and $g_1(\psi(\cq))$ the expressions
\begin{eqnarray*}
\psi(\cq)&=&\cq\prod_{n=1}^\infty(1-\cq^n)^{5(\frac{n}{5})}\\
g_1(\psi(\cq))&=&1+\frac{1}{2}\sum_{n\geq 1} 
((3-i) \chi(n)+(3+i) \overline{\chi(n)}\,) \frac{\cq^n}{1-\cq^n}
\end{eqnarray*}
where $(\frac{.}{5})$ is the quadratic character and $\chi$ is the \emph{quartic character}  
$\left(\frac{.}{2-i}\right)_4$, i.e.
$$
\textstyle{(\frac{n}{5})}=(-1)^j\,,\qquad \chi(n)=i^j\qquad\textrm{if}\quad
n\equiv 2^j\bmod 5;
$$
the latter congruence can also be written as 
$n\equiv i^j\bmod (2-i)$ in $\ZZ[i]$.

From this we infer
$$
\cq\frac{d}{d\cq}\log\CQ=1-\sum_{n\geq 1} b_n \frac{n^2\cq^n}{1-\cq^n}=
1+\frac{1}{2}\sum_{n\geq 1}((2-i) \chi(n)+(2+i) \overline{\chi(n)})
 \frac{n^2\cq^n}{1-\cq^n}\,.
$$

\subsection*{\nr 6: $(A,B,\gl)=(7,-8,2)\,,\;\ga=-1$.}
This example does not appear in \cite{LMW,KMV} or \cite{RV}.
$$
\begin{array}{rrrrr}
n&\qquad&a_n&\qquad&b_n\\
&&&&\\
1&&1&&1\\
2&&-1&&1\\
3&&1&&0\\
4&&-2&&-1\\
5&&5&&-1\\
6&&-14&&0\\
7&&42&&1\\
8&&-136&&1\\
9&& 465&&0\\
\end{array}
$$
Zagier \cite{Z} gives for $\psi(\cq)$ and $g_1(\psi(\cq))$ the expressions
\begin{eqnarray*}
\psi(\cq)&=&\frac{\eta(\cq)^3\eta(\cq^6)^9}{\eta(\cq^2)^3\eta(\cq^3)^9}\\
g_1(\psi(\cq))&=&\frac{\eta(\cq^2)\eta(\cq^3)^6}{\eta(\cq)^2\eta(\cq^6)^3}
=\frac{1}{3}(E_{1,\chi_{-3}}(\cq)+2E_{1,\chi_{-3}}(\cq^2))\,.
\end{eqnarray*}
From this one then obtains
$$
\cq\frac{d}{d\cq}\log\CQ=1-\sum_{n\geq 1} b_n \frac{n^2\cq^n}{1-\cq^n}=
\frac{1}{9}(E_{3,\chi_{-3}}(\cq)+8E_{3,\chi_{-3}}(\cq^2))\,.
$$
Here $E_{1,\chi_{-3}}(\cq)$ and $E_{3,\chi_{-3}}(\cq)$ are the same as in Example \nr 3.\\
Note that 
$b_n= (-1)^{n-1}\chi_{-3}(n)$.

%

%



\begin{thebibliography}{99}
\bibitem{ABYZ}
Ahlgren, S., B. Berndt, A. Yee, A. Zaharescu,
\textit{Integrals of Eisenstein series and derivatives of L-functions},
International Math. Research Notices 32 (2002) 1723--1738.
\bibitem{Beau}
Beauville, A., \textit{Les familles stables de courbes elliptiques sur 
$\mathbf{P}^1$ admettant quatre fibres singuli\`eres}, C. R. Acad. Sc. Paris, t. 294 (1982) 657-660
\bibitem{Bo}
Boyd, D., \textit{Mahler's measure and special values of L-functions},
Experimental Math. vol. 7 (1998) 37--82
\bibitem{CKYZ}
Chiang, T-M., A. Klemm, S-T. Yau, E. Zaslow, \textit{Local mirror symmetry: Calculations and interpretations}, hep-th/9903053
\bibitem{KOS}
Kenyon, R., A. Okounkov, S. Sheffield, \textit{Dimers and Amoebae},
arXiv:math-ph/0311005
\bibitem{KMV}
Klemm, A., P. Mayr, C. Vafa, \textit{BPS states of exceptional non-critical strings}, arXiv:hep-th/9607139
\bibitem{KZ}
Klemm, A., E. Zaslow, \textit{Local Mirror Symmetry at Higher Genus},
arXiv:hep-th/9906046
\bibitem{LMW}
Lerche, W., P. Mayr, N.P. Warner, \textit{Non-critical strings,
Del Pezzo singularities and Seiberg-Witten curves},
arXiv:hep-th/9612085
\bibitem{O}
Okounkov, A., \textit{Random surfaces enumerating algebraic curves},
arXiv:math-ph/0412008
\bibitem{ORV}
Okounkov, A., N. Reshetikin, C. Vafa, \textit{Quantum Calabi-Yau and Classical Crystals},
arXiv:hep-th/0309208
\bibitem{RV}
Rodriguez Villegas, F., \textit{Modular Mahler measures I}, 
 Topics in number theory (University Park, PA, 1997), Ahlgren, S., G. Andrews, K. Ono (eds) 17--48,
Math. Appl., 467,
Kluwer Acad. Publ., Dordrecht, 1999. See 
also: http://www.ma.utexas.edu/users/villegas/research.html
\bibitem{SH}
Schmickler-Hirzebruch, U., \textit{Elliptische Fl\"achen \"uber $\PP_1\CC$ mit drei Ausnahmefasern und die hypergeometrische Differentialgleichung},
Schriftenreihe des Mathematischen Instituts der Universit\"at M\"unster, 2. Serie, Heft 33 (1985)
\bibitem{SB}
Stienstra, J., F. Beukers, \textit{On the Picard-Fuchs equation and the formal Brauer group of certain elliptic K3-surfaces}, Math. Ann. 271 (1985) 269-304
\bibitem{S2}
Stienstra, J., \textit{Mahler Measure, Eisenstein Series and Dimers}, these proceedings
\bibitem{Z}
Zagier, D., \textit{Integral solutions of Ap\'ery-like recurrence equations}, preprint
\end{thebibliography}
\end{document}